\documentclass[11pt, a4]{amsart}
\usepackage[a4paper, left=28mm, right=28mm, top=28mm, bottom=34mm]{geometry}
\usepackage{amsmath} 					%数式環境
\usepackage{amsthm} 						%定理環境
\usepackage{amssymb} 					%記号。txfontsやpxfontsもあり
\usepackage{amscd} 						%簡単な図式。xypicが面倒なときに

\usepackage{epic, eepic}
\usepackage{graphicx} 		%図。上部下部にキャプション可能

\usepackage[all]{xy}							%可換図式ほか。可換図式のラベルに行列を用いるにはxymatrixなどの工夫が要る

\usepackage{mathrsfs} 					%花文字
\usepackage{latexsym}						%証明終わりの記号のための何か
\usepackage{type1ec} % <-- New!
\usepackage[OT2,T1]{fontenc}
\usepackage[russian,english]{babel}
\usepackage{enumerate}					%番号つき箇条書き。descriptionでいい気はするが
\usepackage{verbatim}						%文字通りに書く環境。
\AtBeginDocument{\shorthandoff{"}}			
\newcommand{\Rational}{\mathbb{Q}}					%有理数体
\newcommand{\Real}{\mathbb{R}}					%実数体
\newcommand{\Integer}{\mathbb{Z}}					%整数環
\newcommand{\Complex}{\mathbb{C}}					%複素数体
\def\max{\mathop{\mathrm{max}}\nolimits}			%最大
\newcommand{\isomto}{\overset{\sim}{\to}}

\def\rank{\mathop{\mathrm{rank}}\nolimits}			%階数
\def\det{\mathop{\mathrm{det}}\nolimits}			%行列式
\def\Gal{\mathop{\mathrm{Gal}}\nolimits}		%Galois群
\def\Res{\mathop{\mathrm{Res}}\nolimits}		%制限函手、Weil制限
\def\Cor{\mathop{\mathrm{Cor}}\nolimits}		%余制限準同型
\def\sheafhom{\mathop{\mathscr{H}\kern -2pt om}\nolimits}		%層の射

\def\sheafend{\mathop{\mathscr{E}\kern -2pt nd}\nolimits}			%層の自己準同型
\def\sheafext{\mathop{\mathscr{E}\kern -2pt xt}\nolimits}			%層のExt
\def\Left{\mathop{\mathrm{L} \kern -2pt}\nolimits}				%左導来函手
\def\Right{\mathop{\mathrm{R} \kern -2pt}\nolimits}			%右導来函手
\newcommand{\Cohomology}[2]{H^{#1}\! \left( {#2} \right)}

\newcommand{\strshf}{\mathcal{O}}				%整数環・関数環
\newcommand{\projsp}{\mathbb{P}}				%射影空間。\Pは段落終わり記号
\def\Supp{\mathop{\mathrm{Supp}}\nolimits}			%台
\def\length{\mathop{\mathrm{length}}\nolimits}
\def\Spec{\mathop{\mathrm{Spec}}\nolimits}			%スペクトル
\def\sheafspec{\mathop{\mathscr{S}\kern -2pt pec}\nolimits}
\def\sheafproj{\mathop{\mathscr{P}\kern -2pt roj}\nolimits}

\def\Pic{\mathop{\mathrm{Pic}}\nolimits}			%直線束の同値類
\def\Jac{\mathop{\mathrm{Jac}}\nolimits}		%Jacobians
\def\Brauer{\mathop{\mathrm{Br}}\nolimits}		%Brauer群
\def\GenLin{\mathop{\mathrm{GL}}\nolimits}			%一般線型群
\def\ProjGL{\mathop{\mathrm{PGL}}\nolimits}		%射影線型群
\makeatletter
\def\iddots{\mathinner{\mkern1mu\raise\p@
    \hbox{.}\mkern2mu\raise4\p@\hbox{.}\mkern2mu
    \raise7\p@\vbox{\kern7\p@\hbox{.}}\mkern1mu}}
\def\adots{\mathinner{\mkern2mu\raise\p@\hbox{.} %% yhmath.styから
 \mkern2mu\raise4\p@\hbox{.}\mkern1mu
 \raise7\p@\vbox{\kern7\p@\hbox{.}}\mkern1mu}}
\makeatother

\newtheorem{theorem}{Theorem}[section]
\newtheorem{acknowledgement}[theorem]{Acknowledgement}

\newtheorem{corollary}[theorem] {Corollary}
\newtheorem{definition}[theorem]{Definition}
\newtheorem{example}[theorem]{Example}
\newtheorem{lemma} [theorem]{Lemma}
\newtheorem{notation}[theorem]{Notation}

\newtheorem{proposition}[theorem]{Proposition}
\newtheorem{remark}[theorem]{Remark}

\newtheorem*{prf}{Proof}
\newcommand{\relmiddle}[1]{\mathrel{}\middle#1\mathrel{}}	
\newcommand{\relmid}{\relmiddle{|}}
 \makeatletter
    
    \@addtoreset{equation}{section}
  \makeatother

\newcommand{\Equref}[1]{(\ref{#1})}		%()付きの引用

\title[Proportion of cubic curves with linear determinantal representations]
{A positive proportion of cubic curves over $\mathbb{Q}$ admit
linear determinantal representations}
\date{\today}
\author{Yasuhiro Ishitsuka}
\address{Department of Mathematics, Faculty of Science, Kyoto University, Kyoto 606-8502, Japan}
\email{yasu-ishi@math.kyoto-u.ac.jp}


\begin{document}

\begin{abstract}
Can a smooth plane cubic be defined 
by the determinant of a square matrix with entries in 
linear forms in three variables?  If we can, we say that it admits
a linear determinantal representation. In this paper,
we investigate linear determinantal representations of
smooth plane cubics over various fields, and 
prove that any smooth plane cubic over a large field (or an ample field)
admits a linear determinantal representation.
Since local fields are large, any smooth plane cubic over a local field
always admits a linear determinantal representation.
As an application, we prove that a positive proportion of 
smooth plane cubics over $\mathbb{Q}$, ordered by height,
admit linear determinantal representations. We also prove that, if
the conjecture of Bhargava--Kane--Lenstra--Poonen--Rains
on the distribution of Selmer groups is true, 
a positive proportion of smooth plane cubics over
$\mathbb{Q}$ fail the local-global principle
for the existence of linear determinantal representations.
\end{abstract}

\maketitle

\section{Introduction}
A smooth plane cubic $C \subset \projsp^2$ defined over a field $k$ is
said to \textit{admit a linear determinantal 
representation} over $k$ 
if there is a triple of square matrices $(M_0, M_1, M_2)$
of size three with entries in $k$ such that 
$C$ is defined by the cubic equation
\[
	\det(X_0M_0+X_1M_1+X_2M_2)=0.
\]
It is a classical problem in algebraic geometry to find linear determinantal
representations of plane curves over an algebraically closed field
(\cite[Chapter 4]{Dol12}, \cite{Vin89}).
For concrete examples of conics and cubics, 
see Example \ref{Ex: conic}, Example \ref{Ex: Wei} 
and Example \ref{Ex: Hes}.

Recently, the arithmetic nature of linear 
determinantal representations has been studied by several people 
(\cite{Ho, Ish14, II14a, II14b, II14c}). It is applied to
the computation of the Cassels--Tate pairing of 
elliptic curves (\cite{FN14}).
For an application to the classification of non-associative algebras 
motivated by coding theory, see \cite{DG08, DG11}.

The aim of this paper is to prove that 
a positive proportion of smooth plane cubics over $\Rational$,
ordered by height, admit linear determinantal representations.
We also obtain conditional results on the proportion of smooth plane cubics
over $\Rational$ which fail the local-global principle for the existence of
linear determinantal representations.

\begin{theorem}[See Theorem \ref{Height} and \S 9.4]\label{main4}
	We order all smooth plane cubics over $\Rational$ by height
	(for the definition of height of smooth plane cubics, 
	see Section \ref{CCdet}).
	\begin{enumerate}[(i)]
		\item A positive proportion of smooth plane cubics over $\Rational$ 
		\textit{admit} linear determinantal representations over $\Rational$.
		\item If the conjecture of Bhargava--Kane--Lenstra--Poonen--Rains on
		the distribution of Selmer groups \cite[Conjecture 1.3]{BKLPR15} is true, 
		a positive proportion of smooth plane cubics over $\Rational$ \textit{fail}
		the local-global principle for the existence of 
		linear determinantal representations,
		i.e., they admit linear determinantal representations over each local field
		$\Rational_v$ of $\Rational$, whereas they do not admit
		linear determinantal representations over $\Rational$.
	\end{enumerate}
\end{theorem}

It seems difficult to prove Theorem \ref{main4} (ii) unconditionally.
See Remark \ref{Rmk: Noncond}.

We also obtain unconditional results for linear 
determinantal representations of elliptic curves over $\Rational$.
We regard an elliptic curve $(E, \strshf)$ over $\Rational$ as
a smooth plane cubic by
the complete linear system $\lvert 3 \strshf \rvert \colon E \to \projsp^2$.

\begin{theorem}[See Theorem \ref{main3}]\label{main2}
	We order all isomorphism classes of
	elliptic curves over $\Rational$ by height
	(for the definition of height of elliptic curves over $\Rational$, 
	see Section \ref{Prop}).
	\begin{enumerate}[(i)]
		\item At least 20.68\% of elliptic curves over $\Rational$ 
		\textit{admit} linear determinantal representations over $\Rational$.
		\item At least 16.50\% of elliptic curves over $\Rational$
		\textit{fail}
		the local-global principle for the existence of 
		linear determinantal representations.
	\end{enumerate}
\end{theorem}

These results are based on the global results
recently obtained by Bhargava--Skinner--Zhang \cite{BSZ14},
where they calculate lower bounds of the proportion of elliptic curves over $\Rational$
with Mordell--Weil rank zero (respectively, one).

Also, as a key local result, we prove the existence of linear
determinantal representations over large fields (Theorem \ref{main1}).
Recall that a field $k$ is called \textit{large} (or \textit{ample}) if
any smooth geometrically connected curve over $k$ 
with a $k$-rational point has infinitely many
$k$-rational points (\cite{Pop96, Jar11}). 

\begin{theorem}[{See Theorem \ref{ample}}]\label{main1}
	Any smooth plane cubic over a large field 
	admits a linear determinantal representation.
\end{theorem}

Since local fields (i.e.\ finite extensions of $\Rational_p, \Real$ and
$\mathbb{F}_p((t))$) are large, any smooth plane cubic over a local field
always admits a linear determinantal representation.
Combining the result of Bhargava--Skinner--Zhang, Theorem \ref{main1} and an interpretation of 
linear determinantal representations on smooth plane cubics
using line bundles (see Section \ref{DRandLB}), 
we prove Theorem \ref{main2}. 

Another key result to prove Theorem \ref{main4} is 
the following result on proportion of the smooth plane cubics over $\Rational$.

\begin{theorem}[{See Theorem \ref{Height}}]\label{mainHt}
	We order all smooth plane cubics over $\Rational$ by height
	(for the definition of height of smooth plane cubics, 
	see Section \ref{CCdet}).
	\begin{enumerate}[(i)]
		\item A positive proportion of smooth plane cubics over $\Rational$ have
		Jacobian varieties with positive Mordell--Weil rank.
		\item Assume the conjecture of Bhargava--Kane--Lenstra--Poonen--Rains 
		\cite[Conjecture 1.3]{BKLPR15} is true.
		Then a positive proportion of smooth plane cubics over $\Rational$ have
		Jacobian varieties with trivial Mordell--Weil group.
	\end{enumerate}
\end{theorem}

We prove Theorem \ref{mainHt} using similar methods of Bhargava \cite{Bha14},
where he studies the proportion of smooth plane cubics
over $\Rational$ which satisfy (resp.\ fail) the local-global principle
for the existence of rational points.
We choose a fundamental domain $F$ in the space of ternary cubics over
$\Real$ with respect to the action of $\ProjGL_3(\Integer)$.
Then we prove results similar to Theorem \ref{mainHt}
\textit{with respect to the height of Jacobian varieties
(the Jacobian height)}, not with respect to the height of 
smooth plane cubics (Proposition \ref{JH1}). 
We find that a positive proportion of 
smooth plane cubics in $F$ with integral coefficients 
are \textit{generic} with respect to the Jacobian height.
Since non-generic smooth plane cubics are negligible with respect to
the height of smooth plane cubics,
we can obtain a lower bound of the number of generic plane cubics
with bounded \textit{height} whose Jacobian variety have
positive Mordell--Weil rank 
(resp.\ trivial Mordell--Weil group) 
using the number of generic plane cubics in $F$
with bounded \textit{Jacobian height}
whose Jacobian variety have positive Mordell--Weil rank 
(resp.\ trivial Mordell--Weil group).
Combining Theorem \ref{mainHt} and an interpretation of
linear determinantal representations using line bundles,
we prove Theorem \ref{main4}.

The present paper is organized as follows.
In Section \ref{Prop}, we prove Theorem \ref{main2}
using the results proved in later sections.
In Section \ref{CCdet}, we precisely state Theorem \ref{mainHt}
and prove it using similar methods of \cite{Bha14}.
In Section \ref{Prel}, we recall some basic definitions and facts
about Picard groups and relative Brauer groups.
In Section \ref{DRandLB}, for smooth plane curves of any degree,
we prove a bijection between
equivalence classes of linear determinantal representations and 
isomorphism classes of certain line bundles.
As examples, we study the case of lines and smooth conics 
in Section \ref{LineConic}.
In Section \ref{Cubic}, we obtain some sufficient conditions 
for a smooth plane cubic over a field
to admit a linear determinantal representation.
%We study linear determinantal representations of
%smooth plane cubics over finite fields in Section \ref{FinLDR}.
In Section \ref{LF}, we recall the definition of large fields and
prove Theorem \ref{main1}. Then we apply it to study 
the local-global principle for the existence of 
linear determinantal representations over global fields 
and prove Theorem \ref{main4} in Section \ref{LG}.

\begin{remark}
	These results are in contrast with the case of \textit{symmetric} 
	determinantal representations studied in \cite{II14a} and \cite{II14b}.
	In \cite[Proposition 4.2]{II14a}, it is proved that a smooth plane cubic
	$C$ over a field $k$ admits a symmetric determinantal representation
	if and only if the Jacobian variety $\Jac(C)$ has a non-trivial $k$-rational
	2-torsion point. From this, it is not difficult to prove that
	\begin{itemize}
		\item any smooth plane cubic over a global field \textit{satisfies}
		the local-global principle for the existence of 
		symmetric determinantal representations (\cite[Theorem 5.1]{II14a},
		\cite[Theorem 6.1]{II14b}), and
		\item 100\% of smooth plane cubics over $\Rational$ do \textit{not}
		admit symmetric determinantal representations, and
		\item there are smooth plane cubics over any non-archimedean local field
		without symmetric determinantal representations.
	\end{itemize}
	Hence Theorem \ref{main4}, Theorem \ref{main2} and
	Theorem \ref{main1} do not hold for symmetric determinantal
	representations.	
	See Remark \ref{largeSDR} and Remark \ref{SDR}.
	For a similar problem considering 
	whether a binary form can be written as a discriminant form of 
	a pencil of quadrics, see \cite{BGW13, Creu16}.
\end{remark}

\begin{notation}
For a field $k$, we denote a separable closure of $k$ as $k^s$.
For a $k$-scheme $X$ and a field extension $L$ of $k$,
we denote the base change by $X_L := X \times_{\Spec k} \Spec L$.
%We abbreviate $X^s := X_{k^s}$. 
For a geometrically integral $k$-scheme $X$, 
$k(X)$ denotes the function field of $X$. %and
%$k^s(X)$ denotes that of $X^s$.
For a ring $R$, we denote the multiplicative group as $R^\times$.
\end{notation}

\section{Linear determinantal representations on elliptic curves over $\Rational$}\label{Prop}
In this section, we study the proportion of elliptic curves over $\Rational$ 
which admit (resp.\ do not admit) linear determinantal representations
over $\Rational$.

Before state the precise statement of Theorem \ref{main2}, 
let us recall the notion of height of elliptic curves 
over $\Rational$.
Any elliptic curve $(E, \strshf)$ over $\Rational$ is defined by 
a unique Weierstrass equation
\begin{equation}\label{Ell}
	E_{A,B} \colon (X_1^2X_2 - X_0^3-AX_0X_2^2-BX_2^3 = 0),
\end{equation}
where $A, B \in \Integer$ satisfy $4A^3+27B^2 \neq 0$, and, 
for any prime $p$ with $p^4 \mid A$, 
we have $p^6 \nmid B$. We define
the \textit{height} $H(E_{A,B})$ of the elliptic curve $E_{A,B}$ by
\[
	H(E_{A,B}) := \max \{ 4|A|^3, 27B^2\}.
\]
This is a well-defined invariant for an arbitrary isomorphism class of
elliptic curves over $\Rational$.
%Given two integers $A, B \in \Integer$ satisfying the above conditions,
%we write $E_{A,B}$ defined by \Equref{Ell}.

For a constant $X>0$, we define $\mathrm{Ell}(X)$ 
as the set of isomorphism classes of elliptic curves
over $\Rational$ having height less than $X$.
Let $\mathrm{Ell}(X)_a$ (resp.\ $\mathrm{Ell}(X)_{na}$) be 
the subset of $\mathrm{Ell}(X)$
of elliptic curves over $\Rational$ which admit (resp.\ do not admit) 
linear determinantal representations over $\Rational$.

\begin{theorem}\label{main3}
	We order (isomorphism classes of) elliptic curves 
	over $\Rational$ by height.
	\begin{enumerate}[(i)]
		\item At least 20.68\% of elliptic curves over $\Rational$ 
		\textit{admit}
		linear determinantal representations over $\Rational$, i.e.,
		we have 
		\begin{align*}
			\liminf_{X \to \infty} \frac{\# \mathrm{Ell}(X)_a}
			{\# \mathrm{Ell}(X)} &\ge 0.2068.
		\end{align*}
		\item At least 16.50 \% of elliptic curves over $\Rational$ 
		\textit{do not admit}
		linear determinantal representations over $\Rational$, i.e.,
		we have
		\begin{align*}
			\liminf_{X \to \infty} \frac{\# \mathrm{Ell}(X)_{na}}
			{\# \mathrm{Ell}(X)} &\ge 0.1650.
		\end{align*}
	\end{enumerate}
\end{theorem}

\begin{prf}
(i) 
Bhargava--Skinner--Zhang proved that, when ordered by height,
at least 20.68\% of elliptic curves over $\Rational$ 
have Mordell--Weil rank 1 (\cite[Theorem 3]{BSZ14}).
By Corollary \ref{Jac3}, an elliptic curve over $\Rational$
with positive Mordell--Weil rank admits a linear determinantal
representation over $\Rational$ (actually, it has infinitely many
equivalence classes of 
linear determinantal representations over $\Rational$). 
Hence at least 16.50\% of elliptic curves 
over $\Rational$ admit
linear determinantal representations over $\Rational$.

(ii) Bhargava--Skinner--Zhang proved that, when ordered by height,
at least 16.50\% of elliptic curves over $\Rational$ have
Mordell--Weil rank 0 (\cite[Theorem 3]{BSZ14}). 
It is well-known that 
100\% of elliptic curves over $\Rational$ have Mordell-Weil group
with trivial torsion (for example, see \cite[Lemma 5.7]{BKLPR15}).
Combining these two results, we conclude that
at least 16.50\% of elliptic curves over $\Rational$ have
trivial Mordell--Weil group.
They do not admit linear determinantal representations over $\Rational$
by Lemma \ref{HPforLDR}.
\qed\end{prf}

%\begin{remark}
%	Instead of using the results in \cite{BSh14} and \cite{BSk14},
%	we can use the results of Bhargava--Skinner--W.\ Zhang 
%	in the recent preprint \cite[Theorem 3]{BSZ14} 
%	to improve Theorem \ref{main3}.
%	We can prove that at least 20.68\% of elliptic curves 
%	over $\Rational$ admit linear determinantal representations
%	over $\Rational$, and 	at least 16.50\% of elliptic curves over $\Rational$
%	do not admit linear determinantal representations.
%\end{remark}

\section{Proportion of smooth plane cubics over $\Rational$ and Mordell--Weil group}\label{CCdet}

In this section, we prove that, when ordered by height,
a positive proportion of smooth plane cubics over $\Rational$ have
Jacobian variety with positive Mordell--Weil rank.
We also prove that, if the conjecture of 
Bhargava--Kane--Lenstra--Poonen--Rains 
on the distribution of Selmer groups is true,
a positive proportion of smooth plane cubics over $\Rational$
have Jacobian variety with trivial Mordell--Weil group.
The method of the proof is similar to \cite{Bha14}.

\subsection{Settings and the precise statement of Theorem \ref{mainHt}}
Let $V(\Integer)$ (resp.\ $V(\Rational), V(\Real)$) be the space of ternary cubic forms in three variables 
$X_0, X_1, X_2$ over $\Integer$ (resp.\ $\Rational, \Real$). Since ternary cubic forms can be
parametrized by the coefficients of cubic monomial terms in three variables,
it is isomorphic to the 10-dimensional affine space.
Let us denote the plane cubic over $\Rational$ 
corresponding to $v \in V(\Real)$ as $C(v) \subset \projsp^2_\Real$.

Consider the action of $\GenLin_3(\Real)$ on $V(\Real)$ defined by
\[
	(g \cdot f)(X_0, X_1, X_2) := (\det g)^{-1} f((X_0, X_1, X_2) \cdot g)
\]
for $g \in \GenLin_3(\Real)$ and $v \in V(\Real)$.
This induces an action of $\ProjGL_3(\Real)$ on $V(\Real)$.
It is well-known that the ring of polynomial invariants of 
the representation $(\ProjGL_3(\Real), V(\Real))$ is generated by 
two polynomials $A(v), B(v)$ defined over $\Rational$ of degree 4 and 6, respectively
(see \cite[Theorem 4.4]{Fis08}, \cite[\S 3.2]{AKMMMP}).
Also, there is an invariant 
\[
	\Delta(v) := 4A(v)^3 + 27B(v)^2
\] 
of degree 12 called the \textit{discriminant} of $v$,
which vanishes if and only if the plane cubic $C(v)$ is not smooth.
\begin{remark}\label{Rmk: Sym}
	There are slight differences between our settings 
	and those in \cite{Fis08} and \cite{AKMMMP}.
	In \cite{Fis08}, Fisher considered 
	three invariants $c_4=-A(v)/27$, $c_6=-B(v)/54$ and
	$\Delta = -\Delta(v)/(1728 \cdot 27^3 \cdot 4)$.
	In \cite{AKMMMP}, An--Kim--Marshall--Marshall--McCallum--Pellis 
	considered two invariants $S=A(v)/(1296 \cdot 27)$ and 
	$T=-B(v)/(5832 \cdot 54)$.
\end{remark}

An element $v \in V(\Integer)$ is called \textit{nondegenerate} if $\Delta(v) \neq 0$,
and an element $v \in V(\Integer)$ is called \textit{generic} 
if $v$ is smooth and
$C(v)$ has no $\Rational$-rational flex points (i.e. the $\Rational$-rational points where
the tangent meets the cubic to order at least three).
We note that an elliptic curve $(E, \strshf)$ over $\Rational$
embedded in $\projsp^2$ by the complete linear system $\lvert 3 \strshf
\rvert$ is always non-generic because the origin $\strshf \in E(\Rational)$
is a $\Rational$-rational flex point.
If $v \in V(\Integer)$ is nondegenerate, the Jacobian variety $\Jac(C(v))$
of the smooth plane cubic $C(v)$ is isomorphic to the elliptic curve
\[
	E_{A(v), B(v)} \colon 
	(X_1^2X_2 - X_0^3-A(v)X_0X_2^2-B(v)X_2^3 = 0).
\]

\begin{definition}
	Let $v \in V(\Real)$ be a real ternary cubic form.
	The \textit{height} $H(v)$ of $v$ is the maximum of the absolute values
	of the coefficients of $v$. The \textit{Jacobian height} $H_J(v)$ of $v$ is
	defined as 
	\[
		H_J(v) := \max \{4\lvert A(v)\rvert^3, 27B(v)^2\}.
	\]
\end{definition}
Note that, the Jacobian height of an integral element $v \in V(\Integer)$ 
does not coincide with
the height of the elliptic curve $\Jac(C_v) = E_{A(v), B(v)}$
because we do not reduce the constants $A(v)$ and $B(v)$. 

\begin{definition}
	For a constant $X>0$, we define the subsets
	$V(X), V(X)_{\mathrm{gen}}, V(X)_{\mathrm{rk} \ge 1}, $ 
	$V(X)_{{\mathrm{MW}=0}}$ of $V(\Integer)$ by  
	\begin{align*}
		V(X) &:= \left\{ v \in V(\Integer) \relmid H(v) < X \right\}, \\
		V(X)_{\mathrm{gen}} &:= \left\{ v \in V(X) \relmid v 
		\mbox{ is generic} \right\}, \\
		V(X)_{\mathrm{rk} \ge 1}
		&:= \left\{ v \in V(X)_{\mathrm{gen}} 
		\relmid \rank \Jac(C(v))(\Rational) \ge 1 
		\right\},\\
		V(X)_{\mathrm{MW} = 0} 
		&:= \left\{ v \in V(X)_{\mathrm{gen}} 
		\relmid \Jac(C(v))(\Rational) = 0 
		\right\}.
	\end{align*}
\end{definition}

We fix a fundamental domain $F \subset V(\Real)$ for the
action of $\ProjGL_3(\Integer)$.

\begin{definition}
	For a constant $X>0$, we define the subsets 
	$F_J(X)$, $F_J(X)_{\mathrm{gen}}$, 
	$F_J(X)_{\mathrm{rk} \ge 1}$, $F_J(X)_{\mathrm{MW}=0}$
	\begin{align*}
		F_J(X) &:= 
		\left\{ v \in F \cap V(\Integer) \relmid H_J(v) < X \right\}, \\
		F_J(X)_{\mathrm{gen}} &:= \left\{ v \in F_J(X) \relmid v 
		\mbox{ is generic} \right\}, \\
		F_J(X)_{\mathrm{rk} \ge 1}
		&:= \left\{ v \in F_J(X)_{\mathrm{gen}} 
		\relmid \rank \Jac(C(v))(\Rational) \ge 1 
		\right\},\\
		F_J(X)_{\mathrm{MW} = 0} 
		&:= \left\{ v \in F_J(X)_{\mathrm{gen}} 
		\relmid \Jac(C(v))(\Rational) = 0 
		\right\}.
	\end{align*}
\end{definition}
Now we give the precise statement of the main results of this section.

\begin{theorem}\label{Height}
\begin{enumerate}[(i)]
	\item
	When ordered by height,
	a positive proportion of smooth plane cubics over $\Rational$ 
	have Jacobian varieties with positive Mordell--Weil rank, i.e.,
	we have
	\begin{align*}
		\liminf_{X \to \infty} \frac{\# V(X)_{\mathrm{rk} \ge 1}}
		{\# V(X)} &> 0.
	\end{align*}
	\item
	Assume that the conjecture of Bhargava--Kane--Lenstra--Poonen--Rains
	\cite[Conjecture 1.3]{BKLPR15} is true.
	When ordered by height, 
	a positive proportion of smooth plane cubics over $\Rational$
	have Jacobian varieties with trivial Mordell--Weil group, i.e.,
	we have
	\begin{align*}
		\liminf_{X \to \infty} \frac{\# V(X)_{\mathrm{MW}=0}}
		{\# V(X)} &> 0.
	\end{align*}
\end{enumerate}
\end{theorem}

\subsection{The Bhargava--Kane--Lenstra--Poonen--Rains conjecture}
Before we prove Theorem \ref{Height},
we recall the notion of Selmer groups of elliptic curves 
and some implications of
the conjecture of Bhargava--Kane--Lenstra--Poonen--Rains.

%Let $E$ be an elliptic curve over a field $k$ and fix an integer $n>0$.
%There is a short exact sequence induced by 
%the multiplication-by-$n$ map $[n]_E \colon E \to E$
%\[
%	\xymatrix{
%		0 \ar[r] & E[n](k^s) \ar[r] &
%		E(k^s) \ar[r]^{[n]_E} & E(k^s) \ar[r] & 0.
%	}
%\]
%Taking Galois cohomology, we have an exact sequence
%\[
%	\xymatrix{
%		0 \ar[r] & E(k)/nE(k) \ar[r] &
%		\Cohomology{1}{k, E[n](k^s)} \ar[r] & 
%		\Cohomology{1}{k, E(k^s)}.
%	}
%\]

Let $k$ be a global field, and $n$ an integer invertible in $k$.
The \textit{$n$-Selmer group} $\mathrm{Sel}_n(E)$ of 
an elliptic curve $E$ over $k$ is a finite group 
which parametrizes locally soluble $n$-coverings of $E$
(see \cite[Chapter X]{Sil09} and \cite{CFOSS08}).
%defined as
%the kernel of the composite of maps
%\[
%	\xymatrix{
%		\Cohomology{1}{k, E[n](k^s)} \ar[r]^(0.45){\Res} &
%		\prod_{v} \Cohomology{1}{k_v, E[n](k_v^s)} 
%		\ar[r]
%		&\prod_{v} \Cohomology{1}{k_v, E(k_v^s)},
%	}
%\]
%where $v$ moves all places of $k$, and
%$\Res$ denotes the product of restriction maps for all places $v$.
There exists a short exact sequence of the following form
\begin{equation}\label{Seq: Selmer}
	\xymatrix{
	0 \ar[r] & E(k)\otimes \Integer/n\Integer \ar[r] &
	\mathrm{Sel}_n (E) \ar[r] & \mbox{\cyr {Sh}}(E)[n] \ar[r] & 0,
	}
\end{equation}
where $\mbox{\cyr {Sh}}(E)$ denotes
the \textit{Tate--Shafarevich group} of $E$.
We fix a prime number $p$ invertible in $k$.
Varying $n$ over all powers of $p$, 
we obtain the following short exact sequence
\begin{equation}\label{Seq: Selmer2}
	\xymatrix{
	0 \ar[r] & E(k) \otimes \Rational_p / \Integer_p \ar[r] &
	\mathrm{Sel}_{p^\infty} (E) \ar[r] & 
	\mbox{\cyr {Sh}}(E) [p^\infty] \ar[r] & 0
	}
\end{equation}
from the short exact sequence \Equref{Seq: Selmer}.

Suppose that we are given a short exact sequence $\mathscr{S}$ of
$\Integer_p$-modules. 
The conjecture of Bhargava--Kane--Lenstra--Poonen--Rains
\cite[Conjecture 1.3]{BKLPR15} 
concerns the density of elliptic curves $E$ in $\mathrm{Ell}(X)$
for which the short exact sequence \Equref{Seq: Selmer2} is isomorphic to
$\mathscr{S}$. We do not recall the precise statement of
the conjecture. Instead, we quote from \cite[Section 5.3]{BKLPR15} and 
\cite[Section 5.5]{BKLPR15} the following two implications.

\begin{proposition}\label{Prop: BKLPR}
	Assume that the conjecture of Bhargava--Kane--Lenstra--Poonen--Rains
	\cite[Conjecture 1.3]{BKLPR15} is true.
	\begin{enumerate}[(i)]
		\item When ordered by height, 50\% of elliptic curves over $\Rational$
		have Mordell--Weil rank 0 and 
		50\% have Mordell--Weil rank 1 (\cite[\S 5.3]{BKLPR15}).
		\item When ordered by height, the proportion of 
		elliptic curves over $\Rational$ 	with trivial $p$-Selmer group is
		\[
			\prod_{j \ge 0}(1+p^{-j})^{-1}
		\]
		for any prime number $p$ (\cite[\S 5.5]{BKLPR15}; see also
		\cite[Conjecture 1.1(a)]{PR12}).
	\end{enumerate}
\end{proposition}

We will need the case of Mordell--Weil rank 0 of
Proposition \ref{Prop: BKLPR} (i)
and the case of $p=3$ of Proposition \ref{Prop: BKLPR} (ii).
These results give the following statement for the proportion of
elliptic curves over $\Rational$ with trivial Mordell--Weil group and
non-trivial 3-Selmer group.
Corollary \ref{Cor: 1/8} will be necessary
in the proof of Theorem \ref{main4} (ii).

\begin{corollary}\label{Cor: 1/8}
	Assume that the conjecture of Bhargava--Kane--Lenstra--Poonen--Rains
	\cite[Conjecture 1.3]{BKLPR15} is true.
	When ordered by height, 
	at least $1/8$ of elliptic curves $E$ over $\Rational$ satisfy the
	following conditions:
	\begin{itemize}
		\item the Mordell--Weil group $E(\Rational)$ is trivial, and
		\item the 3-Selmer group $\mathrm{Sel}_3 (E)$ is non-trivial.
	\end{itemize}
\end{corollary}

\begin{prf}
	If the conjecture of
	Bhargava--Kane--Lenstra--Poonen--Rains is true,
	by Proposition \ref{Prop: BKLPR} (i), we have 
	\[
		\lim_{X \to \infty} \frac{\# \left\{ E \in \mathrm{Ell}(X) \relmid
		\rank E(\Rational) = 0 \right\}}{\# \mathrm{Ell}(X)} = \frac{1}{2}.
	\]
	Since 100\% of elliptic curves over $\Rational$ have
	Mordell--Weil group with trivial torsion 
	(cf.\ \cite[Lemma 5.7]{BKLPR15}), we obtain
	\[
		\lim_{X \to \infty} \frac{\# \left\{ E \in \mathrm{Ell}(X) \relmid
		E(\Rational) = 0 \right\}}{\# \mathrm{Ell}(X)} = \frac{1}{2}.
	\]
	On the other hand, by Proposition \ref{Prop: BKLPR} (ii), we have
	\begin{align*}
		\lim_{X \to \infty} \frac{\# \left\{ E \in \mathrm{Ell}(X) \relmid
		\mathrm{Sel}_3(E) = 0 \right\}}{\# \mathrm{Ell}(X)} &= 
		\prod_{j \ge 0} (1+3^{-j})^{-1} \\
		&= 0.319502\dots
		< \frac{3}{8}.
	\end{align*}
	Hence, when ordered by height, 
	at least $1/2-3/8 = 1/8$ of elliptic curves over $\Rational$ have
	trivial Mordell--Weil group and non-trivial 3-Selmer group. 
\qed\end{prf}

\begin{remark}\label{Rmk: Noncond}
	It seems difficult to obtain results like Corollary \ref{Cor: 1/8}
	unconditionally.
	If an elliptic curve $E$ over $\Rational$ has trivial Mordell--Weil group
	and non-trivial 3-Selmer group, its 3-Selmer group $\mathrm{Sel}_3 (E)$
	is expected to have even rank greater than or equal to 2.
	For the moment, it seems difficult to prove that a positive proportion of
	elliptic curves over $\Rational$ have 3-Selmer rank greater than
	or equal to 2. It is the main reason why our Theorem 1.1 (ii)
	is conditional.
\end{remark}

\subsection{The proportion of smooth plane cubics}
First we prove the following proposition concerning 
the proportion of smooth plane cubics over $\Rational$
with respect to the Jacobian height.

\begin{proposition}\label{JH1}
	\begin{enumerate}[(i)]
	\item There is a constant $a >0$ such that
	\begin{align*}
		\# F_J(X)_{\mathrm{rk} \ge 1} &\ge aX^{5/6} + o(X^{5/6}).
	\end{align*}
	\item If the conjecture of Bhargava--Kane--Lenstra--Poonen--Rains 
	\cite[Conjecture 1.3]{BKLPR15} is true, 
	there is a constant $b > 0$ such that
	\begin{align*}
		\# F_J(X)_{\mathrm{MW} = 0} &\ge bX^{5/6} + o(X^{5/6}).
	\end{align*}
	\end{enumerate}
\end{proposition}

\begin{prf}
(i)
	By \cite[Theorem 3.17]{BSh14a}, we see that
%	the number of isomorphism classes of elliptic curves over $\Rational$ 
%	of height less than $X$ is 
	\begin{equation}\label{EllNum}
		\# \mathrm{Ell}(X) = c_1 X^{5/6} + o(X^{5/6})
	\end{equation}
	for some constant $c_1 > 0$.
	Bhargava--Skinner proved in \cite[Theorem 1]{BSk14} that
	\[
		\# \left\{ E \in \mathrm{Ell}(X) \relmid 
		\rank E(\Rational) = 1 \right\} \ge c_2 X^{5/6} + o(X^{5/6})
	\]
	for some constant $c_2 > 0$.
	By the short exact sequence \Equref{Seq: Selmer},
	each elliptic curve $E$ over $\Rational$ with Mordell--Weil rank 1 has
	at least 2 non-trivial elements in the
	3-Selmer group $\mathrm{Sel}_3(E)$.
	Each non-trivial element of the 3-Selmer groups $\mathrm{Sel}_3(E)$ gives 
	at least one generic element $v$ in $F_J(X)_{\mathrm{rk} \ge 1}$.
	(See \cite[Proposition 30]{BSh14}, where generic elements are
	called \textit{strongly irreducible}.)
	Hence the number of generic elements 
	$v \in F_J(X)_{\mathrm{rk} \ge 1}$ is at least
	\[
		2 c_2X^{5/6} + o(X^{5/6}).
	\]
	Putting $a = 2c_2$, we obtain the desired inequality.

	(ii)
	By Corollary \ref{Cor: 1/8}, 
	if the conjecture of Bhargava--Kane--Lenstra--Poonen--Rains is true,
	at least $1/8$ of elliptic curves over $\Rational$ have trivial
	Mordell--Weil group and non-trivial 3-Selmer group.
	Such curve has at least 2 non-trivial elements in the
	3-Selmer group $\mathrm{Sel}_3(E)$, and each non-trivial element of the group 
	$\mathrm{Sel}_3(E)$ gives at least one generic elements in
	$F_J(X)_{\mathrm{MW}=0}$. Thus we obtain
	\begin{align*}\label{EllNum3}
		F_J(X)_{\mathrm{MW}=0} &\ge 2 \cdot \frac{1}{8} \cdot c_1
		 X^{5/6} + o(X^{5/6}) \\
		 &= \frac{c_1}{4} X^{5/6} + o(X^{5/6}),
	\end{align*}
	where $c_1$ is the same constant as in \Equref{EllNum}.
	Putting $b=c_1/4$, we obtain the second inequality.
\qed\end{prf}

We are now ready to prove Theorem \ref{Height}.

\begin{prf}[Proof of Theorem \ref{Height}]
(i)
	Bhargava--Shankar proved in \cite[Theorem 8]{BSh14}
	that
	\begin{equation}\label{Eq: generics}
		\# F_J(X)_{\mathrm{gen}} = c_3 X^{5/6} + o(X^{5/6})
	\end{equation}
	for some constant $c_3 > 0$. Hence, by Proposition \ref{JH1} (i),
	the number of generic elements 
	$v \in F_J(X)_{\mathrm{gen}}$ with 
	$\rank \Jac(C(v))(\Rational) = 0$ is at most
	\begin{equation}\label{EllNum2}
		(c_3-a)X^{5/6} + o(X^{5/6}).
	\end{equation}
	In other words, at most $1-a/c_3$ of elements in
	 $F_J(X)_{\mathrm{gen}}$ 
	have Jacobian varieties	with Mordell--Weil rank 0.

	Next we consider two regions 
	$D(X), D_J(X)$ of $V(\Real) \cong \Real^{10}$
	defined by
	\begin{align*}\label{Regions}
		D(X) &:= \{ v \in V(\Real) \; | \; 
		H(v) < X^{1/12} \}, \\
		D_J(X) &:= \left\{ v \in F \relmid
		H_J(v) < X \right\}.
	\end{align*}
	Note that $D(X)$ is a bounded region.
	On the other hand, the region $D_J(X)$ is not bounded but 
	has finite volume (\cite[\S 2.3]{BSh14}).
	In fact, we have
	\begin{equation}\label{Eq: Num of J}
		\mathrm{vol}(D_J(X)) = c_3 X^{5/6} + o(X^{5/6}),
	\end{equation}
	where $c_3 > 0$ is the same constant as in \Equref{Eq: generics}
	(see \cite[\S 2, (11)]{BSh14} and \cite[Proposition 16]{BSh14};
	see also the paragraph after Theorem 16 in \cite{Bha14}).
	Since the region $D_J(X)$ has finite volume and
	$H_J(Yv) = Y^{12}H_J(v)$ and $H(Yv) = YH(v)$ for $Y > 0$ 
	and $v \in V(\Real)$, 
	we can take a constant $r > 0$ independently of $X$ such that 
	\[
		\frac{\mathrm{vol} (D(rX) \cap D_J(X))}
		{\mathrm{vol} (D_J(X))} > 1-\frac{a}{c_3}.
	\]
	We fix such a constant $r > 0$, and put
	\[
		\delta := \frac{\mathrm{vol} (D(rX) \cap D_J(X))}
		{\mathrm{vol} (D_J(X))}
		- \left(1-\frac{a}{c_3}\right) > 0.
	\]
	The constant $\delta > 0$ does not depend on $X$.
	Note that $F_J(X) = D_J(X) \cap V(\Integer)$ while
	$F(X^{1/12}) = D(X) \cap V(\Integer)$.

	The region $D(rX) \cap D_J(X)$ is bounded and non-generic elements in
	$D(rX) \cap F_J(X)$ are negligible with respect to $X^{5/6}$ 
	(\cite[\S 2.5]{BSh14}). By \Equref{Eq: Num of J}, we have
	\begin{align*}
		\# (D(rX) \cap F_J(X)_{\mathrm{gen}}) 
		&= \mathrm{vol} (D(rX) \cap D_J(X)) + o(X^{5/6})\\
		&= \left(1-\frac{a}{c_3} + \delta\right)\mathrm{vol}(D_J(X))
		 + o(X^{5/6})\\
		&= (c_3 -a + \delta c_3) X^{5/6} + o(X^{5/6}).
	\end{align*}
	In combination with \Equref{EllNum2}, then we obtain 
	\[
		\# (D(rX) \cap F_J(X)_{\mathrm{rk} \ge 1})
		\ge \delta c_3 X^{5/6} + o(X^{5/6}).
	\]
	Replacing $X$ by $r^{-1}X^{12}$, we obtain
	\[
		\# (D(X^{12}) \cap F_J(r^{-1}X^{12})_{\mathrm{rk} \ge 1})
		\ge \frac{\delta c_3}{r^{5/6}} X^{10} + o(X^{10}).
	\]
	Since $V$ is a 10-dimensional affine space, we have
	\begin{align*}
		\# V(X) &=\# (D(X^{12}) \cap V(\Integer)) \\
		&= c_4X^{10} + o(X^{10})
	\end{align*}
	for some constant $c_4>0$.
	Combining these results, we conclude
	\begin{align*}
		\liminf_{X \to \infty} \frac{\# V(X)_{\mathrm{rk} \ge 1}}
		{\# V(X)}
		&\ge
		\liminf_{X \to \infty} 
		\frac{\# (D(X^{12}) \cap F_J(r^{-1}X^{12})_{\mathrm{rk} \ge 1}
		)}{\# V(X)} \\
		&\ge \frac{\delta c_3}{c_4r^{5/6}} > 0.
	\end{align*}
(ii)
	We briefly sketch the proof because the proof of (ii) is similar to the 
	proof of (i). First we count the number of generic elements in 
	$F_J(X)_{\mathrm{gen}}$
	with $\Jac(C(v))(\Rational) \neq 0$.
	Then we take a constant $r'>0$ independently of $X$ so that
	\[
		\frac{\mathrm{vol} (D(r'X) \cap D_J(X))}
		{\mathrm{vol} (D_J(X))} > 1-\frac{b}{c_3}.
	\]
	From this and Proposition \ref{JH1} (ii), we obtain 
	\[
		\# (F_J(r'^{-1}X^{12})_{\mathrm{MW}=0} \cap
		D(X^{12})) \ge
		\frac{\delta' c_3}{r'^{5/6}} X^{10} + o(X^{10})
	\]
	for the constant 
	\[
		\delta' := \frac{\mathrm{vol} (D(r'X) \cap D_J(X))}
		{\mathrm{vol} (D_J(X))} - \left(1-\frac{b}{c_3}\right) > 0.
	\]
	Note that the constant $\delta' > 0$ does not depend on $X$.
	Thus we conclude
	\begin{align*} 
		\liminf_{X \to \infty} \frac{\# V(X)_{\mathrm{MW} = 0}}
		{\# V(X)}
		&\ge \frac{\delta' c_3}{c_4r'^{5/6}} > 0.
	\end{align*}
	This completes the proof of Theorem \ref{Height}.
\qed\end{prf}

\section{Preliminaries on Picard groups and 
relative Brauer groups}\label{Prel}
In this section, we recall definitions and basic properties of Picard groups
and relative Brauer groups, which will be needed in the proof of
the main results of this paper.
For details, see \cite[Chapter 8]{BLR90}, \cite[\S 1]{CK12}
and references therein.

\subsection{Picard functors and Picard schemes}
Let $C$ be a smooth projective geometrically connected curve 
over a field $k$. We denote the \textit{Picard group} of $C$ by
\(
	\Pic(C) = \Cohomology{1}{C, \strshf_C^\times}.
\)
The sheafification of the fppf-presheaf
\[
	(k \mbox{-Schemes})^{\mathrm{op}} 
	\longrightarrow (\mbox{Sets}), \quad T \mapsto \Pic(C_T) 
	= \Pic( C \times_{\Spec k} T).
\]
is called the \textit{relative Picard functor} $\mathcal{P}\mathrm{ic}_{C/k}$
(\cite[\S 8.1]{BLR90}).
This functor is representable by a $k$-group scheme $\Pic_{C/k}$ locally of finite type 
called the \textit{Picard scheme} (\cite[\S 8.2, Theorem 3]{BLR90}).
The Picard scheme $\Pic_{C/k}$ decomposes to the connected components as
\[
	\Pic_{C/k} = \bigsqcup_{d \in \Integer} \Pic^d_{C/k}.
\]
Each components $\Pic^d_{C/k}$ represents the line bundles on $C$
of degree $d$ (\cite[\S 9.3, Theorem 1]{BLR90}).
Especially, its identity component 
\[
	\Jac(C) := \Pic_{C/k}^0 \subset \Pic_{C/k}
\]
is an abelian variety over $k$ called the \textit{Jacobian variety}. 
Its dimension is equal to the genus $g(C)$ of $C$ 
(\cite[\S 9.2, Proposition 3]{BLR90}).
Other components $\Pic^d_{C/k}$ are torsors under $\Jac(C)$.

\subsection{Relative Brauer groups}

For a scheme $X$, let $\Brauer(X)$ denote 
the (absolute) \textit{Brauer group}, 
the group of equivalence classes of Azumaya algebras over $X$.
For an affine scheme $X = \Spec R$, we write
\[
	\Brauer(R) := \Brauer(\Spec R).
\]
If $X$ is a quasi-projective scheme over a field $k$,
the Brauer group $\Brauer(X)$ is identified with
the \textit{cohomological Brauer group}, i.e., the torsion part of 
$H^{2}_{\mathrm{fppf}}\left( X, \mathbb{G}_m \right) =
H^{2}_{\mathrm{\acute{e}t}}\left( X, \mathbb{G}_m \right)$. 
Since schemes appearing in this paper are 
all smooth and projective over a field $k$, we identify these two groups.
We note that, for a global field $k$, the restriction morphisms yield
an injective homomorphism
\[
	\Brauer(k) \to \bigoplus_v \Brauer(k_v),
\]
where $v$ runs over all places of $k$,
by class field theory (\cite[\S 9, \S 10]{Tat65}).

\begin{definition}
	For a morphism of schemes $X \to Y$, 
	the \textit{relative Brauer group} $\Brauer(X/Y)$ is defined by
	the kernel of the pullback homomorphism $\Brauer(Y) \to \Brauer(X)$. 
	We write 
	\[
		\Brauer(X/R) := \Brauer(X/\Spec R).
	\]
\end{definition}

Let $X$ be a smooth projective geometrically connected scheme 
over a field $k$.
By the five-term exact sequence of the Leray spectral sequence,
we obtain the following exact sequence:
\[
	\xymatrix{
		0 \ar[r] & \Pic(X) \ar[r] & \Pic_{X/k}(k) \ar[r]^(0.6){\delta_X} &
		\Brauer(k) \ar[r] & \Brauer(X).
	}
\]
The map $\Brauer(k) \to \Brauer(X)$ is injective 
if $X$ has a $k$-rational point.
By the definition of the relative Brauer group, we have
the following short exact sequence:
\begin{equation}\label{Pic}
	\xymatrix{
		0 \ar[r] & \Pic(X) \ar[r] & \Pic_{X/k}(k) \ar[r]^(.55){\delta_X} &
		\Brauer(X/k) \ar[r] & 0.
	}
\end{equation}

The sequence \Equref{Pic} implies that, for a $k$-rational divisor class
$\alpha \in \Pic_{X/k}(k)$, its image $\delta_X(\alpha)$ 
in $\Brauer(X/k)$ is 
the obstruction to contain a $k$-rational divisor.

\subsection{Exponents of relative Brauer groups}
In this subsection, we treat exponents of relative Brauer groups.
Let $X$ be a smooth projective geometrically connected 
scheme over a field $k$.

\begin{lemma}\label{order-claim}
	Let $L$ be a finite separable extension of $k$, 
	and assume that $X$ has an $L$-rational point.
	Then the extension degree $[L : k]$ annihilates 
	the relative Brauer group $\Brauer(X/k)$;
	in other words, for an element $\alpha \in \Brauer(X/k)$,
	we have $[L:k]\alpha=0$.
\end{lemma}

\begin{prf}
Let us consider the following commutative diagram:
\[
	\xymatrix{
		\Brauer(k) \ar@<2pt>[d]_{\Res_{L/k}} \ar[r]^{i} &
		\Brauer(X) \ar@<2pt>[d]^{\Res_{X_L/X}}\\
		\Brauer(L) \ar@<2pt>[d]_{\Cor_{L/k}} \ar@{^{(}-{>}}[r]^{i_L}
		& \Brauer(X_L) \ar@<2pt>[d]^{\Cor_{X_L/X}} \\
		\Brauer(k) \ar[r]^{i} &
		\Brauer(X).
	}
\]
The horizontal map $i_L$ in the second line is injective
since $X$ has an $L$-rational point.
Take $\alpha \in \Brauer(k)$ and assume $i(\alpha)=0$.
Then we have
\[
	i_L \circ \Res_{L/k}(\alpha) = \Res_{X_L/X} \circ \; i (\alpha) = 0,
\]
so $\Res_{L/k}(\alpha) = 0$ since $i_L$ is injective.
This shows
\begin{align*}
	[L:k]\alpha &= \Cor_{L/k} \circ \Res_{L/k}(\alpha) = 0.
\end{align*}
\\[-30pt]
\qed\end{prf}

%The \textit{index} $\mathrm{ind}(C)$ of $C$ is defined to be 
%the greatest common divisor of extension degrees
%$[k(x) : k]$ where $x \in C$ moves all closed points in $C$.
%If $C$ has genus 1, the index $\mathrm{ind}(C)$ is 
%the greatest common divisor of extension degrees
%$[k(x) : k]$ where $x \in C$ moves all closed points in $C$
%whose residue field $k(x)$ is separable over $k$ (\cite{Lic68}). 
%Hence in this case, 
%by Lemma \ref{order-claim}, we have the following proposition.

%\begin{proposition}\label{Brauer-order} 
%	The index $\mathrm{ind(C)}$ kills each element in $\Brauer(C/k)$;
%	in other words, for an element $\alpha \in \Brauer(C/k)$,
%	we have $\mathrm{ind}(C)\alpha=0$.
%\end{proposition}

As an easy consequence of the above lemma, we obtain 
a sufficient condition under which the relative Brauer group vanishes.

\begin{corollary}\label{BOrat}
	Assume that one of the following conditions is satisfied:
	\begin{itemize}
		\item $\Brauer(k) = 0$, or 
		\item $X$ has a $k$-rational point, or
		\item $k$ is a global field and $X$ has a $k_v$-rational point
		for each place $v$ of $k$. 
	\end{itemize}
	Then we have $\Brauer(X/k) = 0$.
\end{corollary}

For our purposes, we specialize to the case of smooth plane cubics.

\begin{corollary}\label{BOcor}
	Let $C \subset \projsp^2$ be a smooth plane cubic over a field $k$.
	Then each element of $\Brauer(C/k)$ is killed by 3;
	in other words, for an element $\alpha \in \Brauer(C/k)$,
	we have $3\alpha=0$.
\end{corollary}

\begin{prf}
	When $k$ is a finite field, the Brauer group $\Brauer(k)$ vanishes.
	Hence $\Brauer(C/k) = 0$.
	When $k$ is infinite,
	there exists a line $H \subset \projsp^2$ such that
	$C \cap H$ is an \'etale $k$-scheme of order 3 
	by Bertini's theorem \cite{Jou83}.
	Since $C \cap H$ is defined by a separable polynomial of degree 3 over $k$,
	there is a closed point $x \in C \cap H$
	with $[k(x) : k]$ $=$ 1 or 3. 
	Thus each element in $\Brauer(C/k)$ is killed by 3 
	by Lemma \ref{order-claim}.
\qed\end{prf}

\section{Linear determinantal representations and line bundles}\label{DRandLB}

In this section, we recall the definition of 
linear determinantal representations of smooth plane curves and 
their equivalences. Then we interpret the equivalence classes with 
isomorphism classes of certain line bundles.

Let $C \subset \projsp^2$ 
be a smooth plane curve of degree $d$ over a field $k$.
The genus of the curve $C$ is equal to $g(C) = (d-1)(d-2)/2$.
We fix projective coordinates $X_0, X_1, X_2$ of 
the projective plane $\projsp^2$.

\begin{definition}\label{defLDR}
	\begin{itemize}
	\item A \textit{linear determinantal representation} of $C$ is
	a square matrix $M$ of size $d$ with entries in
	$k$-linear forms in $X_0, X_1, X_2$
	such that the equation
	\[
		(\det(M)=0)
	\]
	is a defining equation of $C$.
	
	\item Two linear determinantal representations $M, M'$ of 
	$C$ are said to be \textit{equivalent} if we can write
	\[
		M' = A M B
	\]
	for some $A, B \in \GenLin_{d}(k)$. 
%	We write $M \sim M'$ if
%	two linear determinantal representations $M, M'$ of $C$ are equivalent.
	\end{itemize}
\end{definition}
We say a line bundle on $C$ is \textit{non-effective}
if its only global section is the zero section.
The following theorem gives an interpretation of 
equivalence classes of linear determinantal representations 
in terms of isomorphism classes of non-effective line bundles.
It is well-known at least
when $k$ is an algebraically closed field of characteristic 0
(\cite{Vin89, Bea00, DG11, Ho, Dol12, Ish14}).

\begin{theorem}\label{LDR}
	There is a natural bijection between the following sets:
	\begin{itemize}
		\item the set of equivalence classes of 
		linear determinantal representations of $C$ over $k$, and
		\item the set of isomorphism classes of non-effective line bundles 
		of degree $g(C)-1$ on $C$.
	\end{itemize}
\end{theorem}

\begin{prf}
	We briefly sketch the proof.
	The following argument works over arbitrary fields $k$.
	We follow arguments of Beauville (\cite[Proposition 1.11]{Bea00}); 
	see also \cite[Proposition 2.6]{Ish14}.
	
	Recall that, for a coherent sheaf $\mathcal{F}$ on $\projsp^2$,
	the following two conditions are equivalent 
	(see \cite[Proposition 2.6]{Ish14}).
	\begin{enumerate}
		\item\label{resol1} 
		The sheaf $\mathcal{F}$ has a minimal locally free resolution
		of the form
		\[
			\xymatrix{
				0 \ar[r] & \strshf_{\projsp^2}(-2)^{r} \ar[r]^M &
				\strshf_{\projsp^2}(-1)^{r} \ar[r] & \mathcal{F} \ar[r] & 0
			}
		\]
		for some $r \ge 1$.
		\item\label{resol2} 
		The sheaf $\mathcal{F}$ is arithmetically Cohen--Macaulay,
		and pure of dimension 1.
		It also satisfies
		\[
			\Cohomology{0}{\projsp^2, \mathcal{F}} =
			\Cohomology{1}{\projsp^2, \mathcal{F}}=0.
		\]
	\end{enumerate}
	When these conditions are satisfied, the sheaf $\mathcal{F}$ is 
	the push-forward of 
	a coherent sheaf on the one-dimensional closed subscheme of 
	$\projsp^2$ defined by
	\(
		(\det(M) = 0).
	\)
	
	We go back to the proof of Theorem \ref{LDR}. 
	Let $\pi \colon C \hookrightarrow \projsp^2$ be a smooth plane curve
	of degree $d$, and $\mathcal{L}$ a line bundle on $C$.	
	Then the push-forward 
	\[
		\mathcal{F} := \pi_* \mathcal{L}
	\] 
	is arithmetically Cohen--Macaulay,
	pure of dimension $1$ and $\Supp \mathcal{F} = C$.
	If $\mathcal{L}$ is non-effective of degree $g(C)-1$, we have
	\[
		\Cohomology{0}{\projsp^2, \mathcal{F}} 
		= \Cohomology{1}{\projsp^2, \mathcal{F}} = 0
	\]
	by the Riemann--Roch theorem for $C$.
	Hence $\mathcal{F}$ satisfies \Equref{resol2},
	and $\mathcal{F}$ has a minimal locally free resolution of 
	the form
			\begin{equation}\label{CohToLDR}
				\xymatrix{
					0 \ar[r] & \strshf_{\projsp^2}(-2)^{r} \ar[r]^M &
					\strshf_{\projsp^2}(-1)^{r} \ar[r] & \mathcal{F} \ar[r] & 0
				}
			\end{equation}
	for some $r \ge 1$. Therefore the equation $(\det(M)=0)$ defines $C$
	and we have $r=d$.	
	Hence $M$ gives a linear determinantal representation of $C$ over $k$.
	A straightforward calculation shows that 
	isomorphic line bundles give equivalent $M$'s, and
	non-isomorphic line bundles give non-equivalent $M$'s.
	
	Conversely, a linear determinantal representation $M$ of $C$
	give a coherent sheaf
	$\mathcal{F}$ on $\projsp^2$ by 
	the short exact sequence \Equref{CohToLDR}. 
	Hence $\mathcal{F}$ satisfies the condition \Equref{resol1}, 
	and we can identify $\mathcal{F}$ 
	with the push-forward of a coherent sheaf $\mathcal{L}$ on $C$.
	Since $\mathcal{F}$ is arithmetically Cohen--Macaulay,
	$\mathcal{L}$ is a vector bundle on $C$.
	By computing the Hilbert polynomial of $\mathcal{L}$, we have 
	\[
		\length_{\strshf_{C,\eta}}(\mathcal{L}_\eta)=1
	\]
	for the generic point $\eta$ of $C$.
	Hence $\mathcal{L}$ is a line bundle on $C$.
	By the vanishing of cohomology
	\[
		\Cohomology{0}{C, \mathcal{L}} 
		= \Cohomology{1}{C, \mathcal{L}} = 0,
	\]
	we obtain that $\mathcal{L}$ is non-effective of degree $g(C)-1$
	by the Riemann--Roch theorem.
	It is straightforward to prove that equivalent 
	linear determinantal representations give 
	isomorphic line bundles. This finishes the proof of Theorem \ref{LDR}. 
\qed\end{prf}

\section{Linear determinantal representations 
of lines and conics}\label{LineConic}

As a simple application of Theorem \ref{LDR}, 
let us treat the case of lines and conics.

\subsection{Lines}
Let $C \subset \projsp^2$ be a line defined over a field $k$.
It is defined by a $k$-linear form $l(X_0, X_1, X_2)$.
The linear form $l$ considered as a square matrix of size one
gives a unique equivalence class of
linear determinantal representations of $C$.

Let us reconsider this from the viewpoint of Theorem \ref{LDR}.
We shall find a non-effective line bundle of degree $g(C)-1 = -1$ on $C$.
The degree homomorphism gives an isomorphism 
\[
	\deg \colon \Pic(C) \isomto \Integer \; ; \quad 
	[\mathcal{L}] \mapsto \deg \mathcal{L}.
\]
Let $P \in C(k)$ be a $k$-rational point. 
The line bundle $\strshf_C(-P)$ is non-effective of degree $-1$, so
there is a unique equivalence class of linear determinantal representations
over $k$ corresponding to $\strshf_C(-P)$.

\subsection{Conics}

Next we treat the case of smooth conics.
%We have the following proposition.
Compare the following proposition with \cite[Proposition 4.1]{II14a}, 
where similar results are 
proved for \textit{symmetric} determinantal representations.

\begin{proposition}\label{Con}
	Let $C \subset \projsp^2$ be a smooth plane conic defined over 
	a field $k$. The following conditions are equivalent.
	\begin{enumerate}[(i)]
		\item\label{Con1} The conic $C$ admits a linear determinantal 
		representation over $k$.
		\item\label{Con3} The conic $C$ has a $k$-rational point.
		\item\label{Con4} The conic $C$ is isomorphic to $\projsp^1$ over $k$.
		\item\label{Con5} The relative Brauer group $\Brauer(C/k)$ vanishes.
	\end{enumerate}
	When these conditions are satisfied, 
	the conic $C$ has a unique equivalence class of
	linear determinantal representations over $k$.
\end{proposition}

\begin{prf}
	We first assume \Equref{Con1}. By Theorem \ref{LDR},
	there exists a non-effective line bundle $\mathcal{L}$ on $C$
	of degree $-1$. Then $\mathcal{L}^{-1}$ is a line bundle of degree one.
	The complete linear system $\left\vert \mathcal{L}^{-1} \right\vert$
	gives
	a $k$-isomorphism $C \isomto \projsp^1$.
	This shows \Equref{Con1} $\Rightarrow$ \Equref{Con4}.
	
	The implication \Equref{Con4} $\Rightarrow$ \Equref{Con3} is
	obvious, and \Equref{Con3} $\Rightarrow$ \Equref{Con5} 
	follows from Corollary \ref{BOrat}.
	
	Since the action of the absolute Galois group $\Gal(k^s/k)$ 
	does not change the degree of line bundles, we have
	\begin{align*}
		\Pic_{C/k}(k) &= \Pic_{C/k}(k^s)^{\Gal(k^s/k)} \\
		&= \Pic_{C/k}(k^s) \\
		&\cong \Integer.
	\end{align*}
	In particular, the degree minus one part 
	$\Pic^{-1}_{C/k}(k)$ is a singleton.
	Hence if we assume \Equref{Con5}, we obtain that
	$\Pic^{-1}(C) \cong \Pic^{-1}_{C/k}(k)$ is a singleton and $C$ has
	a unique equivalence class of linear determinantal representations 
	over $k$ by Theorem \ref{LDR}.
	This shows the implication \Equref{Con5} 
	$\Rightarrow$ \Equref{Con1} and the last assertion.
\qed\end{prf}

As seen in the above proof, the relative Brauer group $\Brauer(C/k)$ is
generated by a single element in $\delta_C([\Pic^1_{C/k}(k)])$
where $\delta_C \colon \Pic_{C/k}(k) \to \Brauer(C/k)$ is the map appearing in \eqref{Pic}.
By class field theory, we have the following corollary.
Compare it with \cite[Theorem 5.1 (1)]{II14a}.

\begin{corollary}
	Let $C \subset \projsp^2$ be a smooth plane conic over a global field $k$.
	Then $C$ admits a linear determinantal representation over $k$
	if and only if $C$ admits a linear determinantal representation 
	over each local field $k_v$ of $k$.
\end{corollary}

\begin{prf}
	By Proposition \ref{Con}, $C$ admits a linear determinantal representation
	over $k$ if and only if the relative Brauer group $\Brauer(C/k)$ vanishes.
	By the above observation, $\Brauer(C/k)$ is generated by
	a single element 
	\[
		\alpha_C \in \delta_C([\Pic^1_{C/k}(k)]).
	\]
	If $C$ admits a linear determinantal representation
	over each local field $k_v$ of $k$, the image of $\alpha_C$ 
	in $\Brauer(k_v)$ is zero for each place $v$ of $k$. 
	By class field theory, $\alpha_C$ vanishes in 
	$\Brauer(k)$ (\cite[\S 9, \S 10]{Tat65}). Hence $C$ admits a 
	linear determinantal representation over $k$.
\qed\end{prf}

\begin{example}[{See \cite{Vin89}.}]\label{Ex: conic}
	Let 	
	\(
		C \subset \projsp^2
	\)
	be a smooth conic over a field $k$ with a $k$-rational point $P \in C(k)$.
	By changing the projective coordinates,
	we may assume $P=[0:0:1]$ and 
	\[
		C = (X_0X_2-X_1^2=0).
	\]
	The conic $C$ has a unique equivalence class of 
	linear determinantal representations over $k$ given by the matrix
	\[
		M = 
		\begin{pmatrix}
			X_0 & X_1 \\
			X_1 & X_2
		\end{pmatrix}.
	\]
	Since $M$ is a symmetric matrix,
	we find any linear determinantal representation of $C \subset \projsp^2$ 
	is equivalent to
	a \textit{symmetric} determinantal representation.
\end{example}

%\begin{remark}
%	The equivalence classes of 
%	\textit{symmetric} determinantal representations
%	of $C$ correspond to the line bundles $\mathcal{L}$ on $C$
%	of degree $-1$ satisfying $\mathcal{L} \otimes \mathcal{L} \cong 
%	\omega_C$, where $\omega_C$ is the canonical line bundle of $C$
%	(\cite{Bea00, Ish14}).
%	Any line bundle of degree $-1$ satisfies the condition 
%	since $\omega_C$ is 
%	a unique line bundle on $C$ of degree $-2$ up to isomorphism.
%	Thus we find any linear determinantal representation of $C$ is
%	equivalent to a symmetric determinantal representation.	
%\end{remark}
%

\section{Linear determinantal representations of cubics}\label{Cubic}
In this section, we study linear determinantal representations
of smooth plane cubics over arbitrary fields.
%We immediately obtain the following results from Theorem \ref{LDR}:
\begin{proposition}\label{nonsp}
	Let $C \subset \projsp^2$ be a smooth plane cubic over a field $k$.
	\begin{enumerate}[(i)]
		\item If $C$ admits a linear determinantal representation over $k$, 
		then $\Jac(C)(k) \neq 0$.
		\item If the relative Brauer group $\Brauer(C/k)$ vanishes
		and $\Jac(C)(k) \neq 0$, then
		$C$ admits a linear determinantal representation over $k$.
	\end{enumerate}
\end{proposition}

\begin{prf}
	(i)
	Let $\mathcal{L}$ be a line bundle on $C$
	corresponding to a linear determinantal representation of $C$ over $k$
	by Theorem \ref{LDR}.
	The line bundle $\mathcal{L}$ is non-effective of degree $0$.
	The image $[\mathcal{L}]$ of $\mathcal{L}$ in	$\Jac(C)(k)$
	is a non-trivial $k$-rational point on $\Jac(C)$.\\
	(ii)
	Since $\Brauer(C/k) = 0$, by the short exact sequence \Equref{Pic}, 
	we have an isomorphism 
	\[
		\Pic(C) \overset{\cong}{\longrightarrow} \Pic_{C/k}(k).
	\]
	Hence each non-trivial $k$-rational point of $\Jac(C)$
	corresponds to a non-trivial line bundle $\mathcal{L}$
	on $C$ of degree 0 defined over $k$.	
	The line bundle $\mathcal{L}$ is non-effective,
	hence it gives a linear determinantal representation of $C$
	by Theorem \ref{LDR}.
\qed\end{prf}

\begin{corollary}\label{Cor: BVanish}
	Assume that one of the following conditions is satisfied:
	\begin{itemize}
		\item $\Brauer(k) = 0$, or 
		\item $C$ has a $k$-rational point, or
		\item $k$ is a global field and $X$ has a $k_v$-rational point
		for each place $v$ of $k$.
	\end{itemize}
	Then $C$ admits a linear determinantal representation over $k$
	if and only if $\Jac(C)(k) \neq 0$.
\end{corollary}

\begin{prf}
	The assertion follows from Corollary \ref{BOrat} and 
	Proposition \ref{nonsp} (ii).
\qed\end{prf}

\begin{remark}\label{JacCE}
	If $C$ has no $k$-rational point, 
	there can be a difference between the set of equivalence classes of
	linear determinantal representations over $k$ and 
	the set $\Jac(C)(k) \setminus \{0\}$.
	It turns out that there are distinct smooth plane cubics $C, C'$ 
	over $\Rational$ satisfying the following conditions:
	\begin{itemize} 
		\item $\Jac(C) \cong \Jac(C')$, and
		\item $C$ admits a linear determinantal representation over $\Rational$, 
		and
		\item $C'$ does not admit a linear determinantal representation 
		over $\Rational$.
	\end{itemize}
	For explicit examples, see Example \ref{Teufel}.
\end{remark}

%We give two sufficient conditions 
%to admit a linear determinantal representation.

\begin{corollary}\label{Jac3}
	Let $C \subset \projsp^2$ be a smooth plane cubic over a field $k$. 
	If there is a $k$-rational point $P \in \Jac(C)(k)$
	with $[3]P \neq 0$,
	then the cubic $C$ admits a linear determinantal representation over $k$.
\end{corollary}

\begin{prf}
%	Take a $k$-rational point $P \in \Jac(C)(k)$ with $[3]P \neq 0$.
	We have \[\delta_C([3]P) = 3\delta_C(P) =0\] 
	in the relative Brauer group $\Brauer(C/k)$
	by Corollary \ref{BOcor}.
	Hence $[3]P$ comes from a non-trivial line bundle $\mathcal{L}$ on $C$
	of degree 0 by the short exact sequence \Equref{Pic}.
	The line bundle $\mathcal{L}$ corresponds to 
	a linear determinantal representation of $C$ over $k$ 
	by Theorem \ref{LDR}.
\qed\end{prf}

\begin{corollary}\label{twopts}
	Let $C \subset \projsp^2$ be a smooth plane cubic over a field $k$ 
	with $\#C(k) \ge 2$.
	Then the cubic $C$ admits a linear determinantal representation over $k$.
\end{corollary}

\begin{prf}
	Since $C$ has a $k$-rational point, 
	the existence of a linear determinantal representation of $C$ over $k$ is
	equivalent to $\Jac(C)(k) \neq 0$ by Corollary \ref{Cor: BVanish}. 
	For two distinct $k$-rational points $P \neq Q$ on $C$,
	the line bundle $\strshf_C(P - Q)$ gives 
	a non-zero $k$-rational point on $\Jac(C)$.
\qed\end{prf}

In the rest of this section, we exhibit some examples of
linear determinantal representations of smooth plane cubics.

\begin{example}\label{Ex: Wei}
	Let $k$ be a field of characteristic not equal to 2 nor 3 and 
	\[
		E \colon ( X_1^2X_2 - X_0^3 - aX_0X_2^2 - bX_2^3 = 0)
		\subset \projsp^2
	\]
	an elliptic curve over $k$ with origin $\strshf = [0:1:0]$
	defined by a Weierstrass equation.
	By Theorem \ref{LDR} and Proposition \ref{nonsp},
	there is a bijection between $E(k) \setminus \{\strshf\}$ and
	the set of equivalence classes of linear determinantal representations of
	$E$ over $k$. This can be given explicitly as follows.
	For a $k$-rational point 
	\[
		P = [\lambda: \mu: 1] \in E(k) \setminus \{\strshf\}, 
	\]
	the matrix
	\[
		M_P := \begin{pmatrix}
			X_0 - \lambda X_2 & 0 & -X_1 - \mu X_2 \\
			-X_1 + \mu X_2 & X_0 + \lambda X_2 & (a+\lambda^2) X_2\\
			0 & X_2 & -X_0
		\end{pmatrix}
	\]
	gives a linear determinantal representation of $E$ over $k$.
	Galinat proved in \cite{Gal14} that they give all distinct representatives of
	equivalence classes of linear determinantal representations of $E$ over $k$.
	See \cite{Vin89} for other representatives
	when $k$ is algebraically closed of characteristic not equal to 2 nor 3.
\end{example}

\begin{example}\label{Ex: Hes}
	Let $k$ be a field of characteristic not equal to 2 nor 3 and 
	\[
		C \colon ( X_0^3 + X_1^3 + X_2^3 + \lambda X_0X_1X_2 = 0)
		\subset \projsp^2
	\]
	a smooth plane cubic over $k$ defined by Hesse's normal form.
	Let 
	\[
		P = [a_0 : a_1: a_2] \in C(k)
	\] 
	be a $k$-rational point with $a_0a_1a_2 \neq 0$. 
	Then the Moore matrix
	\[
		M_P := \begin{pmatrix}
			a_0 X_0 & a_1 X_2 & a_2 X_1 \\
			a_1 X_1 & a_2 X_0 & a_0 X_2\\
			a_2 X_2 & a_0 X_1 & a_1 X_0
		\end{pmatrix}
	\]
	gives a linear determinantal representation of $C$ over $k$.
	Buchweitz and Pavlov proved in \cite[Theorem A]{BP15} that
	two Moore matrices $M_P, M_{P'}$ are equivalent if and only if
	two divisors $3P$ and $3P'$ on $C$ are linearly equivalent.
	Moreover, if $k$ is algebraically closed,
	any equivalence class of linear determinantal representations 
	of $C$ over $k$ can be represented by a Moore matrix $M_P$
	corresponding to a point $P = [a_0: a_1 : a_2] \in C(k)$ with 
	$a_0a_1a_2 \neq 0$.
\end{example}

\section{Linear determinantal representations of plane cubics
over large fields}\label{LF}
In this section, we prove that any smooth cubic over a large field
admits a linear determinantal representation.
Let us recall the notion of \textit{large fields} introduced by F.\ Pop.

\begin{definition}[{\cite{Pop96, Jar11}}]\label{largeF}
	A \textit{large field} $k$ is 
	a field having the following property:
	for any smooth algebraic curve $C$ over $k$ with a $k$-rational point, 
	$C$ has infinitely many $k$-rational points.
\end{definition}

Large fields are also called \textit{ample fields} (\cite{Jar11}).
It is known that the class of large fields is quite rich,
and it contains many interesting fields.
For example, local fields, i.e. finite extensions of $\Rational_p, \mathbb{F}_p((t))$  or $\Real$, are
known to be large as well as $k((x, y))= \mathrm{Frac} \; k[[x,y]]$ for any field $k$.
%
%\begin{itemize}
%	\item PAC (Pseudo-Algebraically Closed) fields are large.
%	Hence all separably closed fields are large.
%	\item Fields which are complete with respect to non-trivial absolute values
%	are large. For example, all non-archimedean and archimedean local fields
%	(i.e.\ finite extensions of $\Rational_p, \mathbb{F}_p((t)), \Real$)
%	are large.
%	\item The quotient field of a domain which is Henselian with respect to
%	a non-trivial ideal is large.
%	For example, $k((x,y)) = \mathrm{Frac} \; k[[x,y]]$ is large
%	(\cite{Pop10}).
%	\item A field whose absolute Galois group is a pro-$p$ group
%	for a prime $p$ is large (\cite[Theorem 5.8.3]{Jar11}).
%\end{itemize}

\begin{theorem}\label{ample}
	Let $C \subset \projsp^2$ be a smooth plane cubic over a large field $k$.
	Then $C$ has infinitely many equivalence classes of linear
	determinantal representations over $k$.
\end{theorem}

\begin{prf}
	Since $k$ is large, the Jacobian variety $\Jac(C)$ has 
	infinitely many $k$-rational points.
	Hence $[3]\Jac(C)(k)$ is an infinite set.
	By Corollary \ref{BOcor}, each element in $[3]\Jac(C)(k)$ has 
	the vanishing image under the homomorphism $\delta_C$,
	so any element in $[3]\Jac(C)(k)$ comes from a line bundle
	of degree 0 on $C$ defined over $k$.
	By Corollary \ref{Jac3}, any element $\alpha \in [3]\Jac(C)(k)$
	with $\alpha \neq 0$ 
	gives an equivalence class of 
	linear determinantal representations of $C$ over $k$.
\qed \end{prf}

\begin{corollary}\label{Cor: Locfield}
	Any smooth plane cubic over a local field $k$
	has infinitely many equivalence classes of 
	linear determinantal representations over $k$.
\end{corollary}

When $k$ is a $p$-adic field,
Corollary \ref{Cor: Locfield} was proved by Deajim--Grant
(see \cite[Proposition 4.1]{DG11}).

\begin{remark}
	For a local field $k$, it is possible to prove 
	$\# \Jac(C)(k) = \infty$ directly.
	It is obvious for $k= \Real$ or $\Complex$.
	When $k$ is non-archimedean, it can be proved using formal group laws
	(see \cite[Chapter VII, Proposition 2.2]{Sil09} and
	\cite[Proposition 4.1]{DG11}).
\end{remark}

\begin{remark}\label{notcubic}
	If the degree of a smooth plane curve $C \subset \projsp^2$ is 
	different from 3, Theorem \ref{ample} need not be true; 
	in fact, it does not hold for smooth conics.
	The Severi--Brauer variety of 
	a quaternion division algebra over a local field $k$ is
	a smooth conic without $k$-rational points.
	Hence it does not admit
	a linear determinantal representation by Proposition \ref{Con}.
	For a concrete example, the conic 
	\[
		C=(X_0^2+X_1^2+X_2^2=0)
	\] 
	over $\Real$ does not have an $\Real$-rational point.
	Hence it does not admit a linear determinantal representation over $\Real$
	by Proposition \ref{Con}.
\end{remark}

\begin{remark}\label{largeSDR}
	Smooth plane cubics over a large field need not
	admit \textit{symmetric} determinantal representations.
	In fact, they correspond to
	a special class of line bundles called non-effective 
	\textit{theta characteristics} (see \cite[Proposition 4.2]{Bea00} and
	\cite[Theorem 1.1]{Ish14}).
	They consist a \textit{finite} subscheme of $\Jac(C)$ over $k$.
	Thus the largeness of the base field $k$ does not assure
	the existence of non-effective theta characteristics over $k$.
\end{remark}

\section{Applications to the local-global principle and Proof of Theorem \ref{main4}}\label{LG}

As applications of the results in the previous section, we give
a sufficient condition for a smooth plane cubic over a global field
to fail the local-global principle for
the existence of linear determinantal representations.

\begin{definition}
	A smooth plane cubic $C \subset \projsp^2$ over a global field $k$ is 
	said to \textit{fail the local-global principle for the existence of
	linear determinantal representations} if it satisfies 
	the following two conditions:
	\begin{itemize}
		\item the cubic $C$ \textit{admits} a linear determinantal representation 
		over each local field $k_v$ of $k$, but
		\item the cubic $C$ \textit{does not admit} a linear determinantal
		representation over $k$.
	\end{itemize}
	If $C$ does not fail the local-global principle,
	we say that $C$ \textit{satisfies 
	the local-global principle for the existence of
	linear determinantal representations}.
\end{definition}

Note that any smooth plane cubic $C \subset
\projsp^2$ over a global field $k$ always admits 
a linear determinantal representation over each local field $k_v$
by Corollary \ref{Cor: Locfield}. Hence we have
the following lemma.

\begin{lemma}\label{Lem: Loc}
	A smooth plane cubic $C \subset \projsp^2$ over a global field $k$ 
	satisfies (resp.\ fails) the local-global principle
	for the existence of linear determinantal representations if and only if
	$C$ admits (resp.\ does not admit) 
	a linear determinantal representation over $k$. 
\end{lemma}

\begin{lemma}\label{HPforLDR}
	Let $k$ be a  global field.
	Any smooth plane cubic $C$ over $k$ with $\Jac(C)(k) = 0$ fails
	the local-global principle for the existence of 
	linear determinantal representations.
\end{lemma}

\begin{prf}
	By Proposition \ref{nonsp} (i),
	the cubic $C$ does not admit a 
	linear determinantal representation over $k$.
	Hence it fails the local-global principle for 
	the existence of linear determinantal representations 
	by Lemma \ref{Lem: Loc}.
\qed\end{prf}

\begin{remark}\label{SDR}
	Lemma \ref{HPforLDR} is in contrast with the results in \cite{II14a}
	and \cite{II14b}. In those papers, 
	it is proved that any smooth plane cubic over a global field
	\textit{satisfies} the local-global principle for the existence of 
	\textit{symmetric} determinantal representations
	(see 
	\cite[Theorem 1.1]{II14a} for the case of characteristic not equal to two
	and \cite[Theorem 1.2]{II14b} for the case of characteristic two).
\end{remark}

By Lemma \ref{HPforLDR}, it is easy to construct
examples of smooth plane cubics which fail the local-global principle
as follows. There are infinitely many explicit examples of
elliptic curves over $\Rational$ with trivial Mordell--Weil group; 
for example, see \cite{NH88}. 
Since elliptic curves can be considered as smooth plane cubics,
we can construct infinitely many smooth plane cubics over $\Rational$ 
which fail the local-global principle for the existence of
linear determinantal representations.

%\begin{example}
%	Any smooth plane cubic $C$ over a global field $k$
%	with $\Jac(C)(k) = 0$ admits 
%	no linear determinantal representations over $k$.
%\end{example}

In general, the condition ``$\Jac(C)(k)=0$'' is not a necessary condition for
a smooth plane cubic $C$ over a global field $k$ to fail
the local-global principle for the existence of 
linear determinantal representations.
We conclude this section with examples of 
smooth plane cubics $C$ over $\Rational$
such that $\Jac(C)(\Rational) \neq 0$ but $C$ does not admit a 
linear determinantal representation over $\Rational$.

\begin{example}\label{Teufel}
Let us give infinitely many examples of smooth plane cubics
$C$ over $\Rational$ such that
\begin{itemize}
	\item the cubic $C$ does not admit 
	a linear determinantal representation over $\Rational$, and
	\item the Mordell--Weil group $\Jac(C)(\Rational)$ is non-trivial, and
	\item the relative Brauer group $\Brauer(C/\Rational)$ is non-trivial.
\end{itemize}
Let $\Rational(u)$ be a cyclic extension of $\Rational$ of degree 3
generated by a root $u$ of the equation
\[
	u^3-9u+9=0,
\] 
and $\tau$ a generator of the Galois group $\Gal(\Rational(u)/\Rational)$.
Let us take a prime number $p$
which is inert in $\Rational(u)$. For example, we can take $p=2$ or 5.
There are infinitely many such primes $p$ by Chebotarev's density theorem
\cite[Theorem 13.4]{Neu99}.
Consider the following smooth plane cubic over $\Rational$:
\[
	C_p \colon ( pX_0^3 + p^2X_1^3 - X_2^3=0) \subset \projsp^2.
\]
The Jacobian variety $E_p = \Jac(C_p)$ of $C_p$ has a Weierstrass equation
\[
	E_p \colon \left( X_1^2X_2-X_0^3 + \frac{27}{4}p^6X_2^3=0 \right)
\] 
(\cite[\S 4]{HHW12}), and is isomorphic to the elliptic curve
of Cremona label 27a1 (\cite{Cre})
\[
	E \colon \left( X_1^2X_2 - X_0^3 + \frac{27}{4}X_2^3 = 0 \right)
\]
independently of the choice of $p$.
The Mordell--Weil group $E_p(\Rational)$ is isomorphic to
$\Integer / 3 \Integer$, and generated by
$[6p^2: -9p^3: 2] \in E_p(\Rational)$. 
Let us calculate the relative Brauer group $\Brauer(C_p/\Rational)$ 
of $C_p$.
We apply \cite[Proposition 4.4]{HHW12} with constants 
$a=p, b=p^2,c=-27p^6/4, r=3p^2, s=-9p^3/2$ and 
\begin{align*}
	\xi &= -\frac{9p^3}{2} + \sqrt{\frac{-27p^6}{4}} \\
	&= \frac{-9+3\sqrt{-3}}{2} \cdot p^3.
\end{align*}
Then the relative Brauer group $\Brauer(C_p/\Rational)$ is 
generated by the class
\[
	\delta_{C_p}([6p^2: -9p^3: 2])=[(\Rational(u)/\Rational, \tau, p)] 
	\in \Brauer(\Rational)[3].
\]
Since the prime $p$ is inert in $\Rational(u)$, 
the ideal $p \strshf_{\Rational(u)}$ is a prime ideal,
and the central simple algebra $(\Rational(u)/\Rational, \tau, p)$ is
a division algebra over $\Rational$ ramified at $p$.
Hence the class 
\[
	[(\Rational(u)/\Rational, \tau, p)] \in \Brauer(\Rational)
\] 
is non-trivial. We conclude
\begin{align*}
	\Brauer(C_p/\Rational) &= 
	\langle \delta_{C_p}([6p^2: -9p^3: 2]) \rangle\\
	&= \langle [(\Rational(u)/\Rational, \tau, p)] \rangle \\
	&\cong \Integer / 3\Integer.
\end{align*}
By the short exact sequence \Equref{Pic}, we have
$\Pic^0(C_p) = 0$ and there is no non-effective line bundle 
of degree 0 on $C_p$ defined over $\Rational$.
By Theorem \ref{LDR}, the cubic $C_p$ does not admit a 
linear determinantal representation over $\Rational$.
Because the elements 
\[
	[(\Rational(u)/\Rational, \tau, p)] \in \Brauer(\Rational)
\] 
are distinct for different choices of primes $p$,
the cubics $C_p$ are not isomorphic to each other over $\Rational$.
%As a corollary, $C$ has no $\Rational$-rational points.
%By Theorem \ref{ample}, $C$ admits 
%a linear determinantal representation over 
%local fields $\Real$ and $\Rational_p$.
\end{example}

\begin{remark}
	The smooth plane cubics $C_p$ in Example \ref{Teufel}
	show another phenomenon stated in Remark \ref{JacCE}.
	Smooth plane cubics $C_p$ and the elliptic curve $E$ have
	isomorphic Jacobian varieties, the cubics $C_p$ do not admit 
	linear determinantal representations over $\Rational$.
	But the elliptic curve $E$ admits a linear determinantal representation 
	over $\Rational$ by Corollary \ref{Cor: BVanish}. 
\end{remark}

Finally, we prove Theorem \ref{main4} using Theorem \ref{Height}.

\begin{prf}[Proof of Theorem \ref{main4}]
(i)
	By Theorem \ref{Height} (i), 
	there are a positive proportion of
	smooth plane cubics over $\Rational$ whose Jacobian varieties have
	Mordell--Weil rank 1.
	By Corollary \ref{Jac3}, they admit
	linear determinantal representations over $\Rational$. 

(ii)
	By Theorem \ref{Height} (ii), 
	there are a positive proportion of smooth plane cubics
	over $\Rational$ whose Jacobian varieties have trivial Mordell--Weil group.
	By Lemma \ref{HPforLDR}, they fail the local-global principle
	for the existence of linear determinantal representations.
\qed\end{prf}

\begin{acknowledgement}
The author would like to thank sincerely to Professor Tetsushi Ito 
for various and inspiring comments.
The author is grateful to Professor Takashi Taniguchi for helpful suggestions.
\end{acknowledgement}

\renewcommand{\refname}{References}

\end{document}